\documentclass{IEEEtran}
\usepackage{cite}
\usepackage{amsmath,amssymb,amsfonts,bm}
\usepackage[mathscr]{euscript}
\usepackage{algorithm}
\usepackage{algorithmicx,algpseudocode}
\usepackage{xcolor}

\renewcommand{\d}{\mathrm{d}}
\let\theoremstyle\relax
\usepackage{graphicx}
\let\theoremstyle\relax

\usepackage{amsthm}

\theoremstyle{plain}
\newtheorem{theorem}{Theorem}[section]
\newtheorem{lemma}[theorem]{Lemma}

\theoremstyle{definition}

\theoremstyle{remark}
\newtheorem*{remark}{Remark}

\newcommand{\x}{\mathrm{x}}
\renewcommand{\v}{\mathrm{v}}
\newcommand{\w}{\mathrm{w}}
\renewcommand{\d}{\mathrm{d}}

\makeatletter

\newcommand{\Rmnum}[1]{\expandafter\@slowromancap\romannumeral #1@}

\def\BibTeX{{\rm B\kern-.05em{\sc i\kern-.025em b}\kern-.08em
    T\kern-.1667em\lower.7ex\hbox{E}\kern-.125emX}}

\begin{document}

\title{Distribution Steering for Discrete-Time Uncertain Ensemble Systems}
\author{Guangyu Wu, \IEEEmembership{Member, IEEE}, Panagiotis Tsiotras,  \IEEEmembership{Fellow, IEEE}, Anders Lindquist,  \IEEEmembership{Life Fellow, IEEE}
\thanks{Guangyu Wu was with Department of Automation, Shanghai Jiao Tong University, Shanghai, China (e-mail: guangyu.wu1@outlook.com).}
\thanks{Panagiotis Tsiotras is with the School of Aerospace Engineering and the Institute for Robotics and Intelligent Machines, Georgia Institute of Technology, Atlanta, Georgia (e-mail: tsiotras@gatech.edu).}
\thanks{Anders Lindquist is with School of Artificial Intelligence, Anhui University, Hefei, China (e-mail: alq@kth.se).}}

\maketitle

\begin{abstract}
Ensemble systems appear frequently in many  engineering applications and, as a result, they have become an important research topic in control theory.
These systems are best characterized by the evolution of their underlying state distribution.
Despite the work to date, few results exist dealing with the problem of directly modifying (i.e., ``steering'')
the distribution of an ensemble system.
In addition, in most existing results, the distribution of the states of an ensemble of 
discrete-time systems is assumed to be Gaussian. 
However, in case the system parameters are uncertain, it is not always realistic to assume that the distribution of the system follows a Gaussian distribution, thus complicating the solution of the overall problem.
In this paper, we address the general distribution steering problem for first-order discrete-time ensemble systems, where the distributions of the system parameters 
and the states are arbitrary with finite first few moments. 
Linear system dynamics are considered using the method of power moments to transform the original infinite-dimensional problem into a finite-dimensional one. 
We also propose a control law for the ensuing moment system, which allows us to obtain the power moments of the desired control inputs. 
Finally, we solve the inverse problem to obtain the feasible control inputs from their corresponding power moments. 
We provide a numerical example to validate our theoretical developments.
\end{abstract}


\section{Introduction}
\label{sec:introduction}

This paper addresses the distribution steering problem for first-order discrete-time ensemble stochastic systems. 
In recent years, the necessity to precisely quantify and manage uncertainty in physical systems has spurred a growing interest in the investigation of the evolution of distributions in a stochastic setting~\cite{yu2024optimal, knaup2023covariance,biswal2020stabilization}.
Distribution steering has a long history but it has been garnering increasing attention recently from researchers both from academia and industry owing to its numerous applications in handling uncertainty in a principled manner~\cite{liu2024optimal, liu2022optimal,chen2015optimal,chen2015optimal2,okamoto2018optimal,balci2020covariance}.

The simplest case of distribution control is probably covariance control, 
the earliest research of which can be traced back to a series of articles by Skelton and his students in the late 1980's and early 1990's~\cite{collins1985covariance, collins1987theory, hsieh1990all, xu1992improved}, which explored the assignability of the state covariance via state feedback over an infinite time horizon.
The problem of controlling the covariance over a finite horizon (e.g., the ``steering'' problem) is much more recent.
Among the numerous results in the literature, one should mention
\cite{okamoto2018optimal, okamoto2019optimal, balci2020covariance, liu2022optimal, bakolas2019dynamic} for the discrete-time case or
\cite{chen2015optimal, chen2015optimal2, yongxin2018optimal, liu2022optimal} for the continuous-time case.
In all these works, the initial and terminal distributions are assumed to be Gaussian.
While the Gaussian assumption is generally acceptable when higher-order moments of the distribution are negligible or not important, the same assumption may not be suitable for many other types of distribution steering problems, 
where controlling higher-order moments is critical.
 
One of the most common applications of ensemble systems is in the area of swarm robotics.
Controlling a swarm of robots with a potentially very large number of robots
has become a prominent research topic that has a diverse number of applications, such as environmental monitoring, autonomous construction, warehouse logistics, and precision agriculture, 
to name just a few~\cite{dorigo2021swarm, dorigo2020reflections,rubenstein2013collective,ikumapayi2024swarm,abhang2024swarm}.
In this context, the objective is not to steer a single robot to a specified state but rather to ensure that all agents collectively satisfy certain macroscopic properties.
To achieve this objective, researchers often model the system using a fluid approximation of the multi-agent system, known as the macroscopic or mean-field model. 
By modeling each agent's dynamics as a Markov process, the mean-field behavior of the population is determined through the Liouville equation corresponding to a Markov process~\cite{bakolas2019dynamic, chen2015optimal2}.
Consequently, the system state is the \textit{distribution} of all the agents, which converges weakly to a continuous distribution as the number of agents approaches infinity \cite{zhang2020modeling}. 
Most importantly, the swarm control problem differs from conventional multi-agent control due to the significantly larger number of agents involved. 
Therefore, graph-theoretic approaches, commonly used in multi-agent and networked control problems~\cite{bullo2009distributed, mesbahi2010graph}, do not apply to a group with a large number of agents, due to scalability concerns. 
Designing a control law that scales well with a large group of agents is crucial for the algorithm to be applicable in real-world settings.

Scalability is not the sole concern in swarm robotic applications, however.
Another equally important issue is the absence of Gaussianity in the case of agents with non-trivial dynamics.
Notably, the mean-field approximation results in non-Gaussian
distributions~\cite{elamvazhuthi2019mean, ringh2023mean}.
Representative results of distribution steering problems that do not require a Gaussian assumption can be found in \cite{sivaramakrishnan2022distribution,deshmukh2018mean, elamvazhuthi2018mean, biswal2021decentralized, nodozi2023physics, caluya2021wasserstein}. 
Also, a conventional feedback control strategy,  in the form of a linear function of the system state,
cannot be employed to address the distribution steering problem with general dynamics. 
For instance, as was shown in \cite{elamvazhuthi2018optimal}, the distribution steering problem where the initial and terminal distributions belong to different function classes, cannot be solved using deterministic feedback laws. 
Designing a control law for this type of distribution steering problems, therefore, poses considerable challenges.

Several attempts have been made to address the general distribution steering problem in the literature. 
One such approach, proposed in \cite{sivaramakrishnan2022distribution}, involves using characteristic functions for discrete-time linear systems with general disturbances. 
Another perspective, presented in \cite{biswal2021decentralized}, treats the system evolution purely as a Markov process, where the control inputs serve as transition kernels. 
This allows the control inputs to be selected as random variables, with their density functions representing the transition probabilities, offering greater flexibility in control input design. 
In this paper, we also explore the same concept, and consider the control inputs as transition probabilities, treating them as random variables.

Ensemble control considers problems with intrinsic perturbations in the system parameters \cite{li2009ensemble, li2010ensemble}. 
However, previous results in ensemble control have primarily focused on the controllability of a single agent subject to a specific perturbation. 
Developing a distribution steering scheme for general ensemble systems 
is currently lacking in the literature.

In this paper, we address some of the aforementioned challenges. 
Specifically, we 
propose a control scheme for the general distribution steering problem of a large group of agents, where we only assume the existence and finiteness of the first few moments of the agent state distribution. 
Based on the moment representation of the original infinite-dimensional system, a finite-dimensional reduction is proposed. 
An optimal control scheme via convex optimization is then proposed, yielding the optimal power moments of the control inputs. 
Finally, 
a realization method is revisited to map the power moments of the control inputs to feasible analytic control actions for each agent.
To the best of our knowledge, this is the first attempt to treat the distribution steering problem for general linear ensemble systems.

\section{Distribution steering of linear ensemble systems}
\label{sec:Distribution}

\subsection{Multi-Agent Ensemble System}

We consider an ensemble system consisting of $N$ members (``agents'').
The agent dynamics are linear and are subject to a perturbation in the system parameter. 
The perturbation is independent of the system state.
Since the agents are assumed to be homogeneous, the system dynamics of the $i^\text{th}$ agent take the form
\begin{equation}   \label{type1}
x_{i}(k+1)=a_{i}(k) x_{i}(k)+u_{i}(k),\quad i=1,\ldots,N,
\end{equation}
where $k = 0,1, \ldots, K$ denotes the time step, and
$x_i(k), u_i(k), a_i(k)$ are all scalars.
Let the initial condition be $x_{i}(0) \sim \chi_{0}$, 
for all $i=1,\ldots,N$.
Both the state $x_{i}(k)$ and control $u_{i}(k)$ are random variables with probability distributions $\chi_{k}$ and $\nu_{k}$, respectively, that is, for all $i=1,\ldots,N$,
$x_{i}(k) \sim \chi_{k}$ and $u_{i}(k) \sim \nu_{k}$ for all $k = 0,1,2, \ldots, K$.
The initial state distribution
$\chi_0$ is arbitrary with the first $2n$ 
power moments being finite. 
We also assume that, for each $i=1,\ldots, N$,  $a_{i}(k)$ are random variables independent of $x_{i}(k)$, with realizations drawn from a known common distribution given by $\alpha_{k}$, that is,  $a_{i}(k) \sim \alpha_{k}$.
In the distribution steering problem formulation considered in this work, the identity of each agent is ignored. 
In other words, all agents follow the same statistics, although their particular realizations may differ. 
Moreover, we also assume that the agents are non-interacting and that the size of each agent is negligible, following the standard modeling assumptions of 
mean-field theory~\cite{elamvazhuthi2019mean}.

\subsection{Problem Statement}

We first give the definition of the general distribution steering problem we consider in this paper. 
Provided with an initial probability density function $\chi_{0}$ of $x_{i}(0)$ and a final probability density function $\chi_{f}$ of $x_{i}(K)$, along with the system equation \eqref{type1}, we wish to determine $\nu_{k}$ for $k = 0, \ldots, N-1$, and also a control sequence $\left( u_{i}(0), \ldots, u_{i}(K-1)\right)$ for each agent $i \in\{1,\ldots,N\}$ such that the terminal state distribution satisfies $\chi_K = \chi_{f}$. 
Unlike the conventional distribution steering 
problem~\cite{balci2020covariance}, 
we do not assume that $\chi_{0}$ and $\chi_{f}$ are necessarily Gaussian.
This lack of Gaussianity severely complicates the problem. 
Henceforth, and without a great loss of generality, we assume that $x_{i}(k), a_i(k)$ and $u_{i}(k)$ are defined on the whole $\mathbb{R}$.

Denote, as usual, by $\mathbb{E}\left[\; \cdot \;\right]$ the expectation operator.
The power moments of the state $x_i(k)$ 
up to order $2n$ obey the equation $
\mathbb{E}\left[ x_{i}^{\ell}(k+1) \right]
= \mathbb{E}\left[ \left( a_{i}(k) x_{i}(k)+u_{i}(k) \right)^{\ell} \right] = \sum_{j=0}^{\ell}\binom{\ell}{j}\mathbb{E} [ a_{i}^{j}(k)x_{i}^{j}(k)u_{i}^{\ell-j}(k) ]$.
We note that it is difficult to treat the term $\mathbb{E}[ a_{i}^{j}(k)x_{i}^{j}(k)u_{i}^{\ell-j}(k) ]$ using conventional methods, 
if the control input $u_{i}(k)$ is given as a function of the current state $x_{i}(k)$.
In \cite{wu2023density, wu2022group, wu2023general}, we assumed $a_{i}(k)$ to be identical for each agent $i$ at time $k$ and proposed a scheme where the control inputs are independent of the current system states, so that we can write $\mathbb{E} [ x_{i}^{j}(k)u_{i}^{\ell-j}(k) ] = \mathbb{E} [ x_{i}^{j}(k)  ]\, \mathbb{E} [ u_{i}^{\ell-j}(k) ]$.
Using the algorithm in \cite{wu2023density, wu2022group, wu2023general}, we were able to steer an arbitrary probability distribution to another one, by only assuming the existence of the first few moments for both distributions.
However, these algorithms require that the system dynamics be stable, i.e., $\left|a_{i}(k) \right| < 1$ for all $k =0, \ldots, K$, thus limiting their applicability. 

\subsection{Density Steering of Ensemble Dynamics}

We consider the density steering of the first-order discrete-time ensemble system \eqref{type1} that is not always stable. 
Instead of the open-loop control inputs of the 
schemes presented in \cite{wu2023density, wu2022group, wu2023general}, 
we propose to use a \textit{feedback control law.} 
Specifically, we choose the control input at time step $k$ as 
\begin{equation}  \label{utildeu}
u_{i}(k) = -c(k)a_{i}(k) x_{i}(k) + \tilde{u}_{i}(k),
\end{equation}
where $x_{i}(k)$ and $\tilde{u}_{i}(k)$ are independent random variables, and $c(k)$ is a constant such that $c(k) \in \left[0, 1\right]$. 
As a result, $u_i(k)$ is a random variable selected using \eqref{utildeu}. 
It is a function of the random variables $a_{i}(k), x_{i}(k), \tilde{u}_{i}(k)$. 
Denote the probability distribution of $\tilde{u}(k)$ as $\tilde{\nu}_{k}$. Since $\tilde{u}_{i}(k)$ is not known for each $i^{\text{th}}$ agent, we cannot obtain $u_{i}(k)$ directly from \eqref{utildeu}.
Instead, in order to control the agent group as a whole, we seek to obtain the control distribution $\tilde{\nu}_{k}$ and to sample $\tilde{u}_{i}(k)$ for each $i^{\text{th}}$ agent and each time step $k$.
%
The problem now becomes one of determining $c(0), \ldots, c(K-1)$, and $ \tilde{u}_{i}(0), \ldots, \tilde{u}_{i}(K-1)$ in  \eqref{utildeu} which are generated by the required control distribution $\tilde{\nu}_{k}$.
We emphasize that in the control law \eqref{utildeu} $\tilde{u}_{i}(k)$ is neither a constant nor a function of $x_{i}(k)$.

By denoting $\tilde{a}_{i}(k) = a_{i}(k) - c(k)a_{i}(k)$, the system equation can be written as
\begin{equation} \label{type1new}
x_{i}(k+1)=\tilde{a}_{i}(k)x_{i}(k)+\tilde{u}_{i}(k).
\end{equation}
and \eqref{utildeu} reads $u_{i}(k) = -c(k)/\left(1 - c(k)\right) \tilde{a}_{i}(k) x_{i}(k) + \tilde{u}_{i}(k)$.
During implementation, the control $u_{i}(k)$ is determined by sampling from the distribution of 
the ``group control'' $u(k)$, which is determined by the ``group state'' $x(k)$.
Below we provide the details for deriving the group state $x(k)$ and the group control $u(k)$.

\subsection{Ensemble Group Dynamics}

The control law in \eqref{utildeu} leads to a convenient 
characterization of the group of agents using occupation measures.
The use of occupation measures, allows us to introduce a single equation that conveniently 
characterizes the ``average'' behavior of the ensemble.
For more details, see also~\cite{pitman1977occupation, zhang2020modeling}.

Let $(\Omega, P)$ be the underlying probability space, and let $\boldsymbol{\mathcal{X}} := \left\{ x_{1}, \ldots, x_{N} \right\}$, $\tilde{\boldsymbol{\mathcal{U}}} := \left\{ \tilde{u}_{1}, \ldots, \tilde{u}_{N} \right\}$, $\tilde{\boldsymbol{\mathcal{A}}} := \left\{ \tilde{a}_{1}, \ldots, \tilde{a}_{N} \right\}$. Define the map $\xi: \mathbb{R}^N \times \mathbb{R}^N \times \mathbb{R}^N \rightarrow \mathcal{M}_{+}$ given by $ \xi(\boldsymbol{\mathcal{X}}, \tilde{\boldsymbol{\mathcal{U}}}, \tilde{\boldsymbol{\mathcal{A}}})=\frac{1}{N} \sum_{i=1}^N \delta_{\left(x_i, \tilde{u}_i, \tilde{a}_i\right)}
$, where $\delta_{z}$ denotes the Dirac measure concentrated at point $z$, and $\mathcal{M}_{+}$ is the set of positive Radon measures. Define $\boldsymbol{\mathcal{X}}(k) := \left\{ x_{1}(k), \ldots, x_{N}(k) \right\}$, $\tilde{\boldsymbol{\mathcal{U}}}(k) := \left\{ \tilde{u}_{1}(k), \ldots, \tilde{u}_{N}(k) \right\}$ and $\tilde{\boldsymbol{\mathcal{A}}}(k) := \left\{ \tilde{a}_{1}(k), \ldots, \tilde{a}_{N}(k) \right\}$. 
Then, we have a random measure $\mu_k: \Omega \rightarrow \mathcal{M}_{+}$, given by $\mu_k(\omega)=\xi(\boldsymbol{\mathcal{X}}(k)(\omega), \tilde{\boldsymbol{\mathcal{U}}}(k)(\omega), \tilde{\boldsymbol{\mathcal{A}}}(k)(\omega))$. 
Note that integrating any function $f(\x, \v, \w)$ against the sum of Dirac measures yields $\int f(\x, \v, \w)\, \d \mu_k(\omega)(\x, \v, \w)=\frac{1}{N} \sum_{i=1}^N f\left(x_i(\omega), u_i(\omega), a_i(\omega)\right)$.
For readability, we suppress the dependence on the sample $\omega \in \Omega$ in the following parts of the paper. 
The measure $\mu_k$ induces marginals for the state, control, and system parameter. 
For instance, the state marginal is $\kappa_k = \frac{1}{N} \sum_{i=1}^N \delta_{x_i(k)}$ and describes the distribution of the agents' states at time step $k$. 
Similarly, we define the control and parameter marginals $\beta_k = \frac{1}{N} \sum_{i=1}^N \delta_{\tilde{u}_i(k)}$ and $\eta_k = \frac{1}{N} \sum_{i=1}^N \delta_{\tilde{a}_i(k)}$, respectively.

In order to characterize the aggregate behavior of the agents, 
we need to re-write the system equation \eqref{type1new}, for all $i \in \{ 1, \ldots, N \}$, as a group system equation. 
In this regard, define the random variables corresponding to the measures $\kappa_{k}, \beta_{k}$, and $\eta_{k}$ 
as $x(k), \tilde{u}(k)$ and $\tilde{a}(k)$, respectively; 
namely, let 
$\mathbb{P}\left\{ x(k) = \x \right\} = \kappa_{k}(\x)$, 
$\mathbb{P}\left\{ \tilde{u}(k) = \v \right\} = \beta_{k}(\v)$, and 
$\mathbb{P}\left\{ \tilde{a}(k) = \w \right\} = \eta_{k}(\w)$.
Hence, for any $r \in \mathbb{R}$, we have that
\begin{equation} \label{Pxk1}
\begin{aligned}
 \mathbb{P}\{ x(k+1) &= r \}
=  \kappa_{k+1} (r) = \frac{1}{N} \sum_{i=1}^N \delta_{x_i(k+1)}(r)\\
= & \frac{1}{N} \sum_{i=1}^N \delta_{\tilde{a}_{i}(k)x_{i}(k)+\tilde{u}_{i}(k)}(r),
\end{aligned}
\end{equation}
which is equal to ${1}/{N}$ only when $\tilde{a}_{i}(k)x_{i}(k)+\tilde{u}_{i}(k) = r$, and is zero otherwise.
Similarly,
\begin{equation}
\begin{aligned}
\mathbb{P}\{ & \tilde{a}(k)x(k)  +\tilde{u}(k) = r \}
=  \int_{\mathbb{R}}\delta_{r}(\x  \w+ \v) \, \d\mu_{k}(\x, \v, \w)\\
&= \frac{1}{N} \sum_{i=1}^N\delta_{r}(\x  \w+ \v)\delta_{\left(x_i, u_i, a_i\right)}(\x, \v, \w),
\end{aligned}
\end{equation}
leading to
\begin{equation}
\mathbb{P}\{ \tilde{a}(k)x(k)  +\tilde{u}(k) = r \} = \left\{
\begin{matrix}
{1}/{N}, & r = \tilde{a}_{i}(k)x_{i}(k)+\tilde{u}_{i}(k), \\[4pt] 
0, & \text{otherwise}.
\end{matrix}
\right.
\end{equation}

Therefore, for all $r \in \mathbb{R}$, we have that
$\mathbb{P}\left\{ x(k+1) = r \right\} = \mathbb{P}\left\{ \tilde{a}(k)x(k)+\tilde{u}(k) = r \right\}$. Since the probability values of the random variables $x(k+1)$ and $\tilde{a}(k)x(k)+\tilde{u}(k)$ are the same,
these two random variables obey identical probability laws, 
which leads to the following system equation of the swarm group 
\begin{equation} \label{eom:swarm}
x(k+1)=\tilde{a}(k)x(k)+\tilde{u}(k),
\end{equation}
where the random variables $x(k), \tilde{u}(k)$ and $\tilde{a}(k)$ represent the system state,  the control input, and the system parameter of the swarm group at time step $k$, respectively. 

Note that we do not take into account the situation where a single point in $\mathbb{R}$ is occupied by more than one agent. 
This is because each agent is assumed to have zero volume, making the event of the overlaps of agents to have measure zero. 
Hence, the probability value of each occupation measure at any point on $\mathbb{R}$ cannot exceed $1/N$. 
{We also have $\mathbb{P}\{ \tilde{a}(k) = r \} = \mathbb{P}\{ \left(1 - c(k)\right)a(k) = r \} = 1/N$ when $r = a_{i}(k)$ and is zero otherwise, which leads to
\begin{equation} \label{ak}
    \tilde{a}(k) = \left(1 - c(k)\right)a(k).
\end{equation}
}

\subsection{Moment System}

When steering a large group of agents, controlling directly the state of each agent is challenging, and may be computationally very expensive. 
Instead of controlling each agent individually, we propose to control the moments of the distribution of agents, which encode the macroscopic statistics of the ensemble. 
Using the previous occupation measures, the moments of the system state and the control inputs can be calculated as follows.
\begin{equation}
\begin{aligned}
  \mathbb{E}\left[ x ^{\ell}(k) \right] &
=  \mathbb{E}\left[\int_{\mathbb{R}} \x ^{\ell} \d \kappa_{k}(\x)\right]\\
= & \mathbb{E} \left[ \frac{1}{N}\sum_{i = 1}^{N} x_{i}^{\ell}(k) \right]
= \frac{1}{N}\sum_{i = 1}^{N} \mathbb{E} \left[x_{i}^{\ell}(k) \right]\\
= & \frac{1}{N}\sum_{i = 1}^{N}\int_{\mathbb{R}}x_{i}^{\ell}\chi_{k}(x_{i})\, \d x_{i}
= \int_{\mathbb{R}}x^{\ell}\chi_{k}(x) \, \d x ,
\end{aligned}
\label{xlK}
\end{equation}
where $\ell = 0,1,\ldots$ is a nonnegative integer, such that 
$\ell \leq 2n$. 
The last equality in \eqref{xlK} stems from the fact that $x_{i}(k) \sim \chi_{k}$ for each $i \in \{1, \cdots, N\}$.

\begin{figure*}[h]
\begin{equation}
\begin{aligned}
& \mathbb{E}\left[ x^{j}(k)\tilde{u}^{\ell-j}(k) \right]\\
= & \mathbb{E} \left[ \int_{\mathbb{R} \times \mathbb{R} \times \mathbb{R}} \x^{j} \v^{\ell - j}\d\mu_{k}(\x, \v, \w) \right]
= \mathbb{E}\left[\frac{1}{N}\sum_{i = 1}^{N} x_{i}^{j}(k) \tilde{u}_{i}^{\ell-j}(k) \right]
= \frac{1}{N}\sum_{i = 1}^{N}\mathbb{E}\left[ x_{i}^{j}(k)\right] \mathbb{E} \left[\tilde{u}_{i}^{\ell-j}(k) \right]\\
= & \frac{1}{N} \left( \sum_{i = 1}^{N}\int_{\mathbb{R}}x_{i}^{j}\chi_{k}(x_{i})\, \d x_{i} \int_{\mathbb{R}}\tilde{u}_{i}^{\ell-j}\nu_{k}(\tilde{u}_{i}) \right) \, \d\tilde{u}_{i}
= \int_{\mathbb{R}}x^{j}\chi_{k}(x) \, \d x \int_{\mathbb{R}}\tilde{u}^{\ell-j}\nu_{k}(\tilde{u}) \, \d\tilde{u} = \mathbb{E}\left[ x^{\ell}(k) \right] \mathbb{E}\left[ \tilde{u}^{\ell - j}(k) \right].
\end{aligned}
\label{Exu}
\end{equation}
\hrulefill
\vspace*{4pt}
\end{figure*}

\begin{figure}[htbp]
\centering
\includegraphics[scale=0.2]{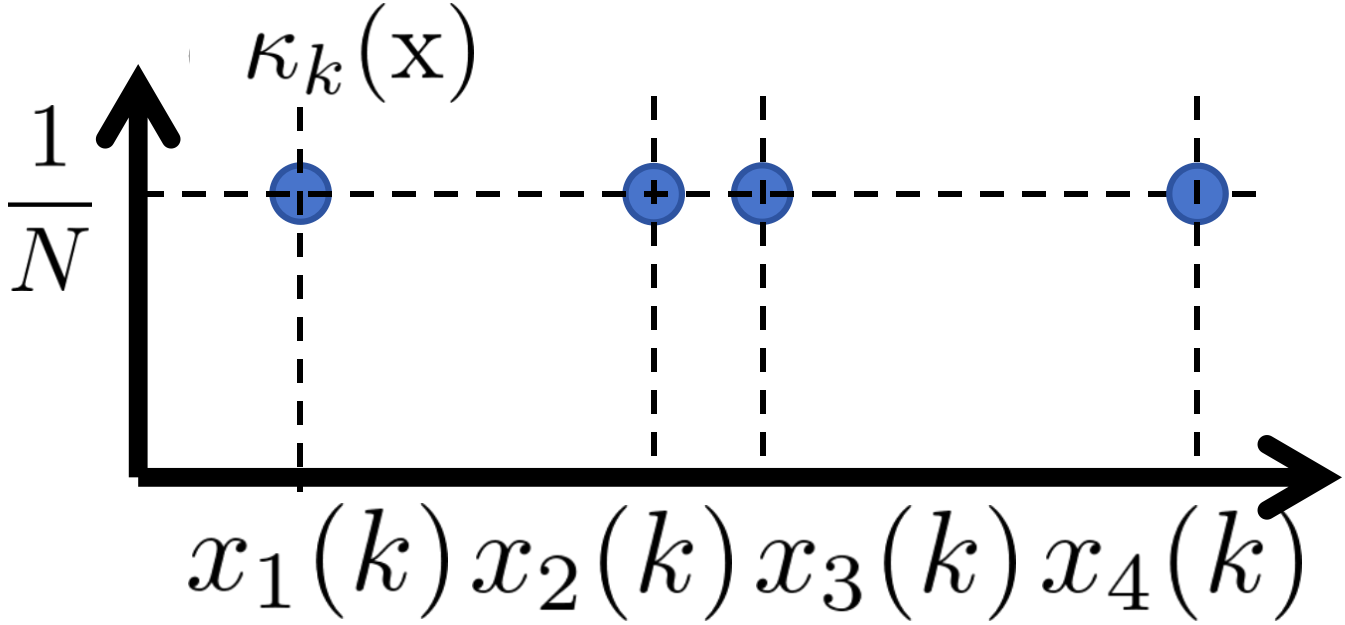}
\centering
\caption{An illustration of the probability distribution $\kappa_{k}(\x)$ of group system state $x(k)$. Each $x_i(k)$ is a sample from $\chi_k(\mathrm{x})$ in $\kappa_k(\mathrm{x})$.}
\label{fig6}
\end{figure}

Then, and similarly to \eqref{xlK}, we obtain
\begin{equation}
\mathbb{E}\left[ \tilde{u}^{\ell}(k) \right]
= \int_{\mathbb{R}}\tilde{u}^{\ell}\tilde{\nu}_{k}(\tilde{u})\, \d \tilde{u},
\label{ulK}
\end{equation}
and $\mathbb{E}\left[ x^{j}(k)\tilde{u}^{\ell-j}(k) \right]$ as in \eqref{Exu}.
The last equality in \eqref{Exu} is owing to the fact $\tilde{u}_{i}(k) \sim \tilde{\nu}_{k}$ for all $i \in \{1, \cdots, N\}$. 
To derive the ensemble dynamics, we take the expectation with respect to the system parameters, effectively integrating $a(k)$ out when calculating the power moments.

Since each $a_{i}(k)$ is independent of $\tilde{u}_{i}(k)$ and $x_{i}(k)$, by following a similar treatment as in \eqref{Exu}, we have that $a(k)$ is independent of $\tilde{u}(k)$ and $x(k)$.
Using \eqref{utildeu} and \eqref{Exu}, we have
\begin{equation}
\mathbb{E}\left[u^{\ell}(k)\right]
=  \sum_{i = 0}^{\ell}\left(-c(k)\right)^{i}\mathbb{E}\left[a^{i}(k)\right]\mathbb{E}\left[x^{i}(k)\right] \mathbb{E}\left[\tilde{u}^{\ell - i}(k)\right].
\label{Euk}
\end{equation}

Then, the dynamics of the moments can be written as the linear matrix equation 
\begin{equation}
    \mathscr{X}(k+1) = \tilde{\mathscr{A}}(\tilde{\mathscr{U}}(k))\mathscr{X}(k)+\tilde{\mathscr{U}}(k),
\label{momentsystem1}
\end{equation}
which is the moment counterpart of the original system \eqref{type1}, where the
system matrix $\tilde{\mathscr{A}}(\tilde{\mathscr{U}}(k))$  is given in \eqref{longeq1}.
Accordingly,
the state vector of \eqref{momentsystem1} is given by
\begin{equation}
\mathscr{X}(k) = \begin{bmatrix}
\mathbb{E}[x(k)] & \mathbb{E}[x^{2}(k)] & \cdots & \mathbb{E}[x^{2n}(k)]
\end{bmatrix}^{\intercal},
\label{XK}
\end{equation}
and the control vector is given by
\begin{equation}
\tilde{\mathscr{U}}(k) = \begin{bmatrix}
\mathbb{E}[\tilde{u}(k)] & \mathbb{E}[\tilde{u}^{2}(k)] & \cdots & \mathbb{E}[\tilde{u}^{2n}(k)]
\end{bmatrix}^{\intercal}.
\end{equation}

\begin{figure*}[t]
\begin{equation}
\tilde{\mathscr{A}}(\tilde{\mathscr{U}}(k))
= \begin{bmatrix}
\mathbb{E}\left[\tilde{a}(k)\right] & 0 & 0 & \cdots & 0\\ 
2\mathbb{E}\left[\tilde{a}(k)\right]\mathbb{E}[\tilde{u}(k)] & \mathbb{E}\left[\tilde{a}^{2}(k)\right] & 0 & \cdots & 0\\ 
3\mathbb{E}\left[\tilde{a}(k)\right]\mathbb{E}[\tilde{u}^{2}(k)] & 3\mathbb{E}\left[\tilde{a}^{2}(k)\right]\mathbb{E}[\tilde{u}(k)] & \mathbb{E}\left[\tilde{a}^{3}(k)\right] & \cdots & 0\\ 
\vdots & \vdots & \vdots & \ddots & \vdots \\ 
\binom{2n}{1}\mathbb{E}\left[\tilde{a}(k)\right]\mathbb{E}[\tilde{u}^{2n-1}(k)] & \binom{2n}{2}\mathbb{E}\left[\tilde{a}^{2}(k)\right]\mathbb{E}[\tilde{u}^{2n-2}(k)] & \binom{2n}{3}\mathbb{E}\left[\tilde{a}^{3}(k)\right]\mathbb{E}[\tilde{u}^{2n-3}(k)] &  \cdots & \mathbb{E}\left[\tilde{a}^{2n}(k)\right]
\end{bmatrix}
\label{longeq1}
\end{equation}
\hrulefill
\vspace*{4pt}
\end{figure*}

 We note that the form of the moment system \eqref{momentsystem1} is similar to the one we proposed in our previous work \cite{wu2022group}. The only difference is that the parameters in \eqref{longeq1} are the power moments of $\tilde{a}(k)$. 
 We note that only finitely many orders of power moments appear in \eqref{longeq1}. 
 Hence, using the moment representation \eqref{momentsystem1}, the original problem which is infinite-dimensional, 
 is approximated with a finite-dimensional one.
 Using \eqref{momentsystem1}, the original distribution steering problem is reduced to a problem of steering the corresponding moment system, which is formulated as follows.

 \subsection{Moment Steering Problem}
 
 Given an arbitrary initial density $\chi_{0}$, determine the control sequences $c(0), \ldots, c(K-1) $ and $u(0), \ldots, u(K-1)$ such that,
for all $i=1,\ldots,N$,
the first $2n$ order moments of the final density $\chi_K$ of $x_i(K)$ are identical to the moments of a specified final probability density, that is, for $\ell = 1, \cdots, 2n$, $\int_{\mathbb{R}} x^{\ell} \chi_{K}(x) \, \d x = \int_{\mathbb{R}} x^{\ell} \chi_{f}(x) \, \d x$.

In this paper, we consider maximizing the smoothness of the state transition in terms of the power moments \cite{wu2023general},
which leads to the following optimization problem
\begin{equation}
\min_{\mathscr{X}(1), \ldots, \mathscr{X}(K-1)} \mathcal{L} \left( \mathscr{X}(1), \ldots, \mathscr{X}(K-1) \right).
\label{OptiProb}
\end{equation}
where,  
$
\mathcal{L} \left( \mathscr{X}(1), \ldots, \mathscr{X}(K-1) \right) := \sum_{k = 0}^{K-1}\big( \mathscr{X}(k+1)  - \mathscr{X}(k) \big)^{\intercal}$
$\big( \mathscr{X}(k+1) - \mathscr{X}(k) \big).
$
The directional derivative of $\delta \mathcal{L} $ along $\delta \mathscr{X}(k)$ then reads
$$
\begin{aligned}
 \delta \mathcal{L}  \big( \mathscr{X}(1), &\ldots, \mathscr{X}(K-1); \delta \mathscr{X}(k) \big)\\
= & 2 \left( \mathscr{X}(k) - \mathscr{X}(k-1) \right) - 2 \left( \mathscr{X}(k+1) - \mathscr{X}(k) \right).
\end{aligned}
$$
The directional derivative of $\delta \mathcal{L}$ has to be zero at a minimum for all variations $\delta \mathscr{X}(k)$, for each $k \in \left\{ 0, \cdots, K-1 \right\}$. 
Therefore, for all $k = 0, \ldots, K-1$, we have $
\mathscr{X}(k) - \mathscr{X}(k-1) = \mathscr{X}(k+1) - \mathscr{X}(k)$. 
It is easy to verify that
\begin{equation}
\mathscr{X}(k) = \frac{K - k}{K}\, \mathscr{X}(0) + \frac{k}{K} \mathscr{X}(K).
\label{Smooth}
\end{equation}
It is important to note that \eqref{OptiProb} is not the only approach for determining the trajectory of the moments; however, it has demonstrated promising performance in various distribution steering tasks. 
Further research will focus on refining metrics to enhance its performance. 
The power moments of the system states of the original system \eqref{type1} up to order $2n$ are then determined for all $k=1,\ldots,K-1$ from \eqref{Smooth}.
However, the existence of $x(k)$ for $k = 1, \ldots, K-1$ given the moments in \eqref{Smooth} still remains to be shown. Next, we provide a proof of the existence of such an $x(k)$. 

\begin{lemma}
Given the moment sequence $\mathscr{X}(0), \ldots, \mathscr{X}(K-1)$ satisfying \eqref{Smooth}, there always exists a state sequence of the original system $x(0), \ldots, x(K-1)$, not necessarily unique, which corresponds to this moment sequence.
\label{Lemma32}
\end{lemma}

\begin{proof}
The statement is equivalent to proving that, for all $\mathscr{X}(k)$ satisfying \eqref{Smooth}, there exists $x(k)$ satisfying \eqref{XK} for all $k = 1, \ldots, K-1$. 
To this end, first, define the Hankel matrix 
\begin{equation}
[\mathscr{X}(k)]_{H} := \begin{bmatrix}
1 & \mathbb{E}[x(k)] & \cdots & \mathbb{E}[x^{n}(k)]\\ 
\mathbb{E}[x(k)] &  \mathbb{E}[x^{2}(k)]& \cdots & \mathbb{E}[x^{n+1}(k)]\\ 
\vdots & \vdots & \ddots & \\ 
\mathbb{E}[x^{n}(k)] & \mathbb{E}[x^{n+1}(k)] &  & \mathbb{E}[x^{2n}(k)]
\end{bmatrix}.
\end{equation}
From~\cite[Theorem~3.8]{schmudgen2017moment}, it suffices to prove that, for all $1 \leq k \leq K-1$,
$[\mathscr{X}(k)]_{H} \succ 0$.
Since $\chi_{0}, \chi_{K}$ are specified initial and terminal densities, $\mathscr{X}(0), \mathscr{X}(K)$ exist. 
It follows that $[\mathscr{X}(0)]_{H}, [\mathscr{X}(K)]_{H}$ are both positive definite.
Using \eqref{Smooth}, by rearranging the elements of $\mathscr{X}(k)$, we have
\begin{equation} \label{XkSumHankel}
[\mathscr{X}(k)]_{H} = \frac{K - k}{K} \, [\mathscr{X}(0)]_{H} + \frac{k}{K} [\mathscr{X}(K)]_{H}.
\end{equation}
Since the scalars $\left(K - k\right)/K$ and  $k/K$ are both positive for all $k = 1, \ldots, K-1$, it follows that $[\mathscr{X}(k)]_{H}$, being the sum of two positive definite matrices, is also positive definite, thus completing the proof. 
\end{proof}

\begin{remark}
We should also note that 
\begin{equation}     \label{qkxinf}
\chi_{k}(x) = \frac{K - k}{K} \chi_{0}(x) + \frac{k}{K} \chi_{K}(x),
\end{equation}
is a feasible distribution of $x(k)$.
To see this, and by denoting the $\ell^{\rm th}$ element of $\mathscr{X}(k)$ as $\mathscr{X}_{\ell}(k)$ we have, for $\ell = 1, \ldots, 2n$,  that $
\mathscr{X}_{\ell}(k) = \int_{\mathbb{R}} x^{\ell} \chi_{k}(x) \, \d x
= \int_{\mathbb{R}}x^{\ell }\,  \frac{K - k}{K} \chi_{0}(x)\, \d x + \int_{\mathbb{R}} x^{\ell} \, \frac{k}{K} \chi_{K}(x) \,\d x
= \frac{K - k}{K} \, \mathscr{X}_{\ell}(0) + \frac{k}{K}  \mathscr{X}_{\ell}(K)$,
which leads to \eqref{XkSumHankel}.
Therefore, the distribution in \eqref{qkxinf} satisfies the moment condition \eqref{Smooth}.
However, there is an infinite number of feasible $x(k)$ corresponding to a given $\mathscr{X}(k)$. 
Furthermore, for distribution steering tasks, we desire the distribution of the system state at each time step to be analytic and the parameter space to be 
finite-dimensional~\cite{glass2009infinite}, which yields a feasible control input $u(k)$.
Finding such an analytic $\chi_{k}$ for each $[\mathscr{X}(k)]_{H} \succ 0$ is treated in the next section.
\end{remark}

\section{Optimal Controller Design}    \label{sec3}

\subsection{Existence of Density Steering Controller}

We have proposed a novel feedback control law, given in \eqref{utildeu}, for solving the distribution steering task. 
We note that designing such a control is considerably easier than the treatments in \cite{wu2023density, wu2022group, wu2023general}, since $c(k) = 1$ is always a feasible solution. 
By setting $c(k) = 1$  for all $ k = 0, \ldots, K-1$, it is always possible to obtain a feasible control input $u(k)$ given $\mathscr{X}(k)$. 
Based on this observation, we propose an optimization scheme for determining the control sequences 
$\mathbf{c} = ( c(0), \ldots, c(K-1)) $ and $\boldsymbol{\tilde{\mathscr{U}}} = \big( \tilde{\mathscr{U}}(0), \ldots, \tilde{\mathscr{U}}(K-1) \big)$.

\begin{theorem}    \label{theorem34}
The optimization problem 
\begin{subequations}   \label{optimization1}
\begin{align}
& \min_{0 \leq c(k) \leq 1} \mathbb{E}\big[\left( -c(k)a(k)x(k) + \tilde{u}(k) \right)^{2} \big],  \label{optimization1A}\\
\mathrm{s.t.} \ & \tilde{\mathscr{U}}(k) = \mathscr{X}(k+1) - \tilde{\mathscr{A}}(\tilde{\mathscr{U}}(k))\mathscr{X}(k), \ [\tilde{\mathscr{U}}(k)]_{H} \succ 0,
\label{optimization1B}
\end{align}
\end{subequations}
is convex.
\end{theorem}

\begin{proof}
We first prove that the cost function \eqref{optimization1A} is convex. 
The second-order power moment of $u(k)$ yields $\mathbb{E}\left[ u^{2}(k) \right] = c^{2}(k)\mathbb{E}\left[a^{2}(k)\right]\mathbb{E}\left[x^{2}(k)\right] - 2c(k)\mathbb{E}\left[a(k)\right]\mathbb{E}\left[x(k)\right]\mathbb{E}\left[\tilde{u}(k)\right] + \mathbb{E}\left[\tilde{u}^{2}(k)\right]$.
Noting that $\frac{\mathrm{d}^{2}\mathbb{E}\left[ u^{2}(k) \right]}{\mathrm{d}c(k)^{2}} = \mathbb{E}\left[a^{2}(k)\right]\mathbb{E}\left[ x^{2}(k) \right] \geq 0$, yields that the cost function is convex~\cite{boyd2004convex}. 
Next, we prove that the domain is a convex set. We need to prove that the feasible domain of $c(k)$ under the constraints in \eqref{optimization1} is convex. 
First, note that $[\tilde{\mathscr{U}}(k)]_{H}$ is a continuous matrix function of $c(k)$. 
Hence, there exists $0< \epsilon < 1 $ such that, for all
$c(k) \in ( \epsilon, 1]$,
$[\tilde{\mathscr{U}}(k)]_{H} \succ 0$.
However, there may be several subintervals of 
$[0, 1]$ that are not path-connected that 
satisfy $[ \tilde{\mathscr{U}}(k) ]_{H} \succ 0$. 
Under this circumstance, the domain of $c(k)$ is not convex. Therefore, we need to prove that there exists an $\epsilon > 0$, such that $[ \tilde{\mathscr{U}}(k)]_{H} \succ 0$ for $c(k) > \epsilon$, and $[ \tilde{\mathscr{U}}(k) ]_{H} \nsucc 0$ for $c(k) < \epsilon$. 
This is equivalent to showing that there exists a feasible $\tilde{u}(k)$ for $c(k) > \epsilon$, while there is no feasible $\tilde{u}(k)$ for $c(k) < \epsilon$.

Assume that $\tilde{u}_{1}(k), \tilde{u}_{2}(k)$ exist given $c_{1}(k), c_{2}(k)$. 
We need to show that, for any $c_{3}(k) \in \left[ c_{1}(k), c_{2}(k) \right]$, $\tilde{u}_{3}(k)$ exists. 
First, write $c_{3}(k) = \lambda c_{1}(k) + (1 - \lambda) c_{2}(k), 0 \leq \lambda \leq 1$.
From equations \eqref{eom:swarm} and \eqref{ak}, 
we have $x(k+1)
= \left( 1 -c_{1}(k) \right)a(k)x(k) + \tilde{u}_{1}(k)
= \left( 1 -c_{2}(k) \right)a(k)x(k) + \tilde{u}_{2}(k)$.

Therefore, we can write
\begin{equation}
\begin{aligned}
x(k+1)
= & \lambda \left( \left( 1 -c_{1}(k) \right)a(k)x(k) + \tilde{u}_{1}(k)\right)\\
+ & (1 - \lambda) \left(\left( 1 -c_{2}(k) \right)a(k)x(k) + \tilde{u}_{2}(k)\right)\\
= & \left( 1 - \lambda c_{1}(k) - \left( 1 - \lambda \right)c_{2}(k) \right)a(k)x(k)\\
+ & \lambda \tilde{u}_{1}(k) + \left( 1 - \lambda \right)\tilde{u}_{2}(k)\\
= & \left( 1 - c_{3}(k) \right)a(k)x(k) + \lambda \tilde{u}_{1}(k)\\
+ & \left( 1 - \lambda \right)\tilde{u}_{2}(k).
\end{aligned}
\label{xk1lambda}
\end{equation}
It follows that 
$\tilde{u}_3(k) = \lambda \tilde{u}_{1}(k) + \left( 1 - \lambda \right)\tilde{u}_{2}(k)$ is the control corresponding to $c_{3}(k)$ and  hence 
the set of all feasible $c(k)$ is convex. 
We also conclude that there exists an $\epsilon > 0$ such that for any $c(k) \in \left( \epsilon, 1 \right]$, $u(k)$ exists, i.e., $[ \tilde{\mathscr{U}}(k)]_{H} \succ 0$. 
\end{proof}

It remains to show the existence of a solution to the optimization problem \eqref{optimization1}. 
In case the conditions of \eqref{optimization1B} are satisfied with $c(k) = 0$,
the feasible set of $c(k)$ is $[0, 1]$, which is closed and convex, 
which ensures the existence of a solution to the optimization problem. 
Otherwise, let the feasible domain of $c(k)$ be $\left(\epsilon, 1\right]$, where $0 < \epsilon < 1$. 
In this case,  we may relax the second condition in \eqref{optimization1B} to $[\tilde{\mathscr{U}}(k)]_{H} \succeq 0$. 
Since $[\tilde{\mathscr{U}}(k)]_{H}$ is a continuous matrix function of $c(k)$, the feasible domain of $c(k)$ is the closed and convex set $[ \epsilon, 1 ]$. 
Hence, the existence of a solution to the optimization problem follows. 
Furthermore, from \eqref{xk1lambda}, $[\tilde{\mathscr{U}}(k)]_{H}$ can be singular only at $c(k) = \epsilon$. 
When $c(k) = \epsilon$ is the optimal solution to the optimization problem \eqref{optimization1}, $\tilde{\nu}(k)$ is an atomic distribution supported on $n$ discrete points on $\mathbb{R}$, rather than a continuous distribution.

Theorem~\ref{theorem34} along with the proof of the existence of a solution allows us to obtain an optimal control $\tilde{\mathscr{U}}(k)$ for each $k = 0, \ldots, K-1$.
However, $\tilde{\mathscr{U}}(k)$ consists of the statistics of the random variable $\tilde{u}(k)$. 
The problem now becomes one of determining a control $\tilde{u}(k)$ given $\tilde{\mathscr{U}}(k)$ obtained by the solution to the optimization problem~\eqref{optimization1}.
In our previous work~\cite{wu2023density, wu2022group, wu2023general}, this step is called the \textit{realization problem} of the random variables $\tilde{u}(k)$ for $k = 0, \ldots, K-1$. 
Here we adopt the same treatment found in \cite{wu2022group}.

For the sake of simplicity, henceforth, we omit the index $k$ if there is no danger of ambiguity. 
The problem now becomes one of designing an algorithm to estimate the probability density supported on $\mathbb{R}$ of which the power moments are given. 
This is known in the literature as the
Hamburger moment problem~\cite{schmudgen2017moment}. 
Often, the Kullback-Leibler divergence is a widely used measure in the literature 
to characterize the difference between the reference density and the density estimate~\cite{byrnes2001finite, byrnes2002identifiability, georgiou2003kullback}. 
A convex optimization scheme for density estimation using the Kullback-Leibler divergence has been proposed in \cite{wu2023non} for the Hamburger moment problem. 
We adopt this strategy to realize the control inputs. 

Let $\mathcal{P}$ be the space of probability density functions defined and having support on the real line, and let $\mathcal{P}_{2n}$ be the subset of all $p \in \mathcal{P}$ that have at least $2n$ finite moments.
The Kullback-Leibler divergence between the probability density functions $p,r \in \mathcal{P}$ is defined as $\mathsf{KL} (r \| p) := \int_{\mathbb{R}} r(\tilde{u}) \log \frac{r(\tilde{u})}{p(\tilde{u})} \, \d \tilde{u}$.
Define the linear operator $\Gamma : \mathcal{P}_{2n} \rightarrow \mathbb{R}^{(n+1) \times (n+1)}$ as $\Gamma ( \tilde{p} ) = \Sigma := \int_{\mathbb{R}} G(\tilde{u}) \tilde{p}(\tilde{u}) G^{\intercal}(\tilde{u}) \mathrm{d}\tilde{u}$, where
$
G(\tilde{u})= \begin{bmatrix}
1 & \tilde{u} & \cdots & \tilde{u}^{n}
\end{bmatrix}^{\intercal}$.
It can be easily shown that 
\begin{equation}  \label{eq:SIGMA}
\Sigma = \begin{bmatrix}
1 & \mathbb{E}[\tilde{u}] & \cdots & \mathbb{E}[\tilde{u}]\\ 
\mathbb{E}[\tilde{u}] &  \mathbb{E}[\tilde{u}^{2}]& \cdots & \mathbb{E}[\tilde{u}^{n+1}]\\ 
\vdots & \vdots & \ddots & \\ 
\mathbb{E}[\tilde{u}^{n}] & \mathbb{E}[\tilde{u}^{n+1}] &  & \mathbb{E}[\tilde{u}^{2n}]
\end{bmatrix},
\end{equation}
where $\mathbb{E}[\tilde{u}^{\ell}],~(\ell = 1, \ldots, 2n)$ are obtained 
from \eqref{momentsystem1} using $\mathscr{X}(1), \ldots, \mathscr{X}(K-1)$ obtained by the optimization problem \eqref{OptiProb}.

Since $\mathcal{P}_{2n}$ is convex, $\operatorname{range}(\Gamma)=\Gamma\mathcal{P}_{2n}$ is also convex.
Given $r \in \mathcal{P}$ and $\Sigma \succ 0$, then there is a unique $\tilde{\nu}^* \in \mathcal{P}_{2 n}$ that minimizes the Kullback-Leibler distance $\mathsf{KL} (r \| \tilde{\nu})$ 
subject to $\Gamma ( \tilde{\nu} ) =\Sigma$, which is given by
\begin{equation}  \label{hatp}
\tilde{\nu}^*=\frac{r}{G^{\intercal} {\Lambda}^* G}, 
\end{equation}
where ${\Lambda}^*$ is the unique solution to the minimization problem~\cite{wu2023non}
\begin{equation}   \label{Jr}
\hspace*{-2mm}
\min_{ \Lambda \in \mathcal{L}_{+}   } \mathcal{J}_{r}(\Lambda):=\operatorname{tr}(\Lambda \Sigma)-\int_{\mathbb{R}} r(\tilde{u}) \log \left[G(\tilde{u})^{\intercal} \Lambda G(\tilde{u})\right] \mathrm{d}\tilde{u},
\end{equation}
where $\mathcal{L}_{+}:=\left\{\Lambda \in \operatorname{range}(\Gamma) \mid G(\tilde{u})^{\intercal} \Lambda G(\tilde{u})>0, \tilde{u} \in \mathbb{R}\right\}$.

The probability density function of the random variable $\tilde{u}$ can now be estimated by solving the convex optimization problem~\eqref{Jr}. 
From Theorem \ref{theorem34}, 
we obtain the values of $c(k)$ for all  $k = 0, \ldots, K-1$ by solving the convex optimization problem \eqref{optimization1}.
Therefore, the control input $u(k)$ can be uniquely determined by solving the two convex optimization problems \eqref{optimization1} and \eqref{Jr}. 
It is worth noting that the power moments of the proposed density estimate align exactly with the specified moments. 
This property distinguishes the proposed approach from other similar moment-matching methods in the literature~\cite{matyas1999generalized}.
Consequently, the proposed approach can be used to realize the control inputs. 
Since both the prior density $r(\tilde{u})$ and the density $\tilde{\nu}^*(\tilde{u})$ are supported on the real line, one can usually select a Gaussian distribution for $r(\tilde{u})$ when $\chi_{f}$ is a sub-Gaussian distribution \cite{wu2023non}, or a Cauchy distribution when $\chi_{f}$ is heavy-tailed.

\subsection{Algorithm for Density Steering}

We now propose an algorithm for steering a large, but finite, group of agents. 
The distribution of the system state is discrete, representing the individual agents. 
We do not aim to steer a specific agent to a specific state. 
Instead, we target the terminal discrete system state distribution to be the desired one. 

It should be noted that the control law is not a purely state feedback control law.
Instead, each control input is the sum of a state feedback function and a random variable that is independent of the current system state. 
A similar idea has appeared in \cite{biswal2021decentralized}, where the evolution of the state distribution is regarded as a Markov process with the control input serving as the transition rate or probability. 

An algorithm for solving the discrete distribution steering is given in Algorithm~\ref{alg:2}. 
Since the random part of the control input $\tilde{u}(k)$ and the current system state $x(k)$ are independent, we may obtain each $\tilde{u}_{i}(k)$ by drawing i.i.d samples from the realized distribution 
$\tilde{\nu}^*_{k}(\tilde{u})$. 
By doing this, $\tilde{u}(k)$ serves as a transition probability of a Markov process.

\begin{algorithm}
    \caption{Discrete distribution steering of a large group of agents}
    \label{alg:2}
    \begin{algorithmic}[1]
        \Require Number of agents $N$; maximal time step $K$; system parameters $a_{i}(k) \sim \alpha_{k}$ for each $i^{\text{th}}$ agent at $k = 0, \ldots, K-1$; 
        initial discrete distribution $\chi_{0}$; 
        specified terminal discrete distribution $\chi_{f}$
        \Ensure Control inputs of the $i^\text{th}$ agent $u_{i}(k)$, $k = 0, \ldots, K-1$, $i = 1, \ldots, N$
        \State $k \leftarrow 0$
        \State Calculate $\mathscr{X}(0)$ from \eqref{XK}
    \While{$0 \leq k < K$}
        \State Calculate $\mathscr{X}(k+1)$ from \eqref{Smooth}.
        \State Solve optimization problem~\eqref{optimization1} with $\mathscr{X}(k+1)$, $a(k)$, $\mathbb{E}[x(k)]$ and $\mathbb{E}[x^{2}(k)]$ (namely, the first two elements of $\mathscr{X}(k)$), to obtain $c(k)$ and $\tilde{\mathscr{U}}(k)$ and $\Sigma$ from \eqref{eq:SIGMA}
        \State Optimize cost function \eqref{Jr} and obtain $\tilde{\nu}^*_{k}$ from~\eqref{hatp}
        \State Draw $N$ i.i.d. samples $\tilde{u}_{i}(k)$ for $i = 1, \ldots, N$, from the distribution $\tilde{\nu}^*_{k}$
        \State Calculate control inputs $u_{i}(k)$ for each agent $i = 1, \ldots, N$ from \eqref{utildeu}
        \State $k \leftarrow k+1$.
    \EndWhile
    \end{algorithmic}
\end{algorithm}

\section{A numerical example}
\label{sec5}

In this section, we provide a numerical example to demonstrate the proposed approach.
We first simulate the discrete-time Liouville control problem, where an infinite number of agents is assumed, and the control inputs $u(k)$ are continuous functions. 
However, since $x(0)$ and $u(k)$ for all $k=0,1,\ldots, K$ are not assumed to fall within an exponential family (e.g., Gaussian), the probability density function of $x(k)$ for $k = 1, \ldots, K$ does not always have an analytic form. 
This makes comparing the results using our algorithm to the desired distribution a difficult task. 
To validate the performance of the proposed algorithm, we simulated a large group of $5000$ agents, based on the results of the Liouville control problem. 
The solution is a sufficiently fine discrete distribution, which makes it possible to compare it to the desired continuous distribution. 
We simulate a steering problem in four steps ($K=4$) where the system parameter follows a Laplace distribution. 
The initial density is chosen as 
$
\chi_{0}(x) = \frac{1}{\sqrt{2\pi}}e^{\frac{x ^{2}}{2}}.
$
The terminal density function is specified as a multi-modal density 
which is a mixture of two Gaussian densities as follows
\begin{equation}
     \chi_{f}(x) = \frac{0.5}{\sqrt{2\pi}}e^{\frac{(x + 2) ^{2}}{2}} + \frac{0.5}{\sqrt{2\pi}}e^{\frac{(x - 3) ^{2}}{2}}.
\label{qt2}
\end{equation}

\begin{figure}[htbp]
\centering
\includegraphics[scale=0.3]{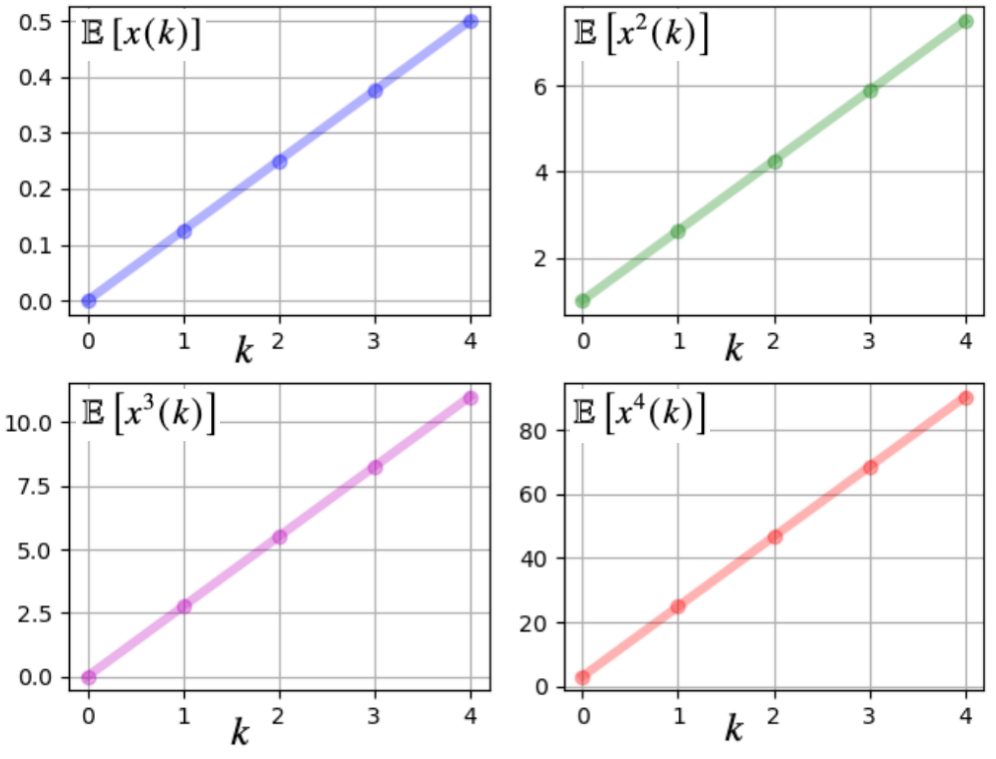}
\centering
\caption{$\mathscr{X}(k)$ at time steps $k = 0, 1, 2, 3, 4$.}
\label{fig1}
\end{figure}

\begin{figure}[htbp]
\centering
\includegraphics[scale=0.3]{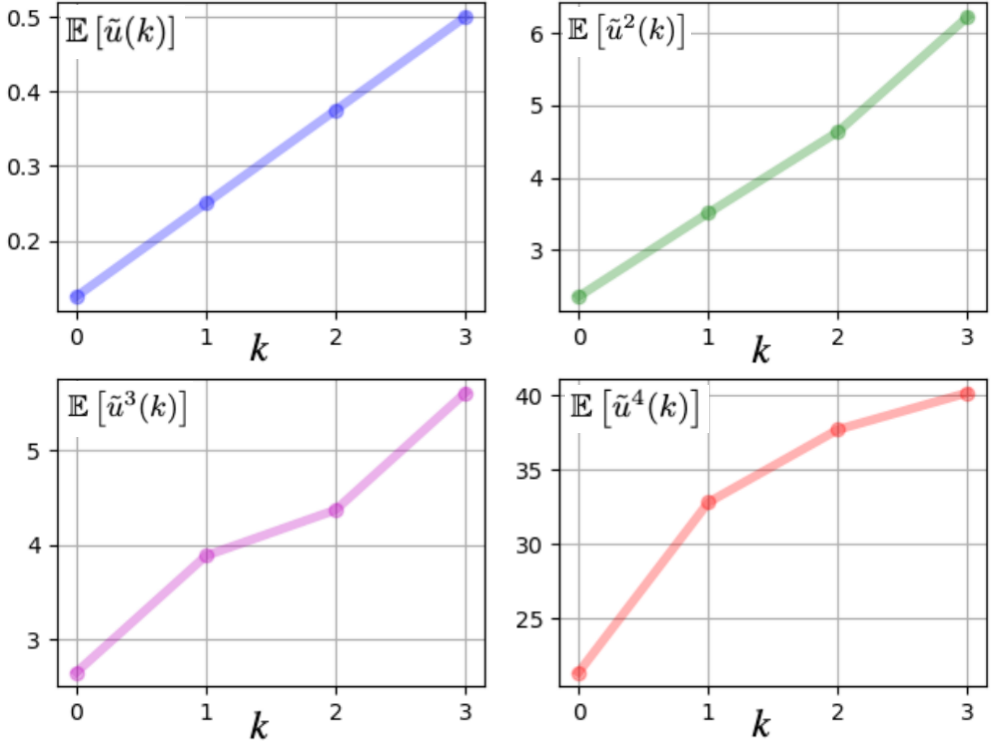}
\centering
\caption{$\tilde{\mathscr{U}}(k)$ at time steps $k = 0, 1, 2, 3$.}
\label{fig2}
\end{figure}

\begin{figure}[htbp]
\centering
\includegraphics[scale=0.3]{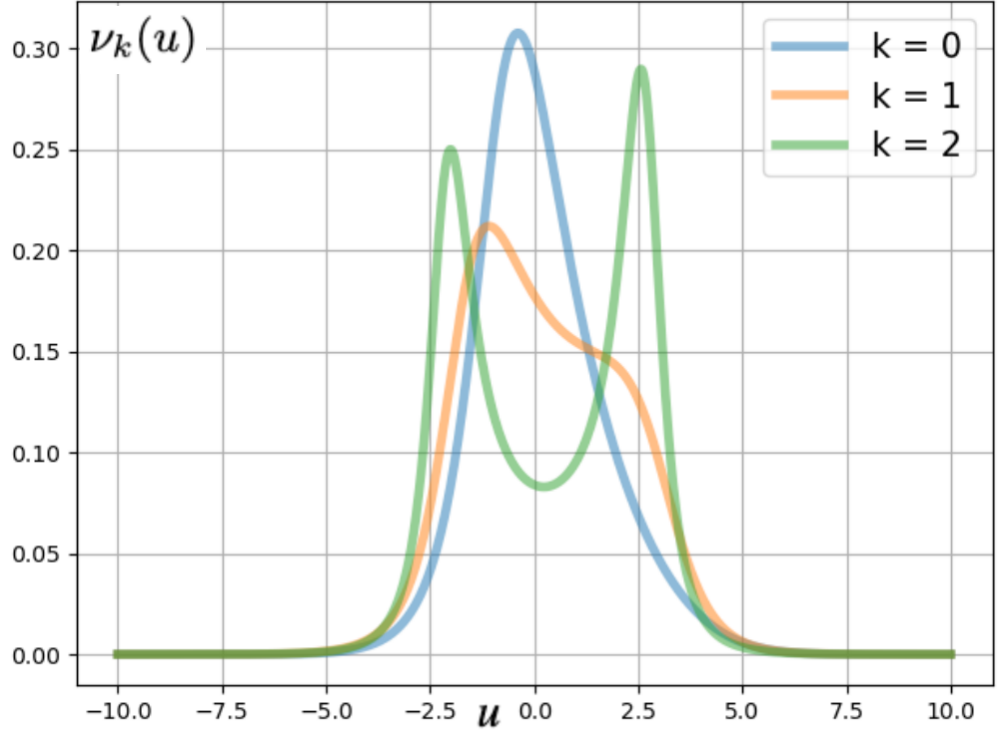}
\centering
\caption{Realized distributions of control inputs $\nu_{k}(u)$ by $\tilde{\mathscr{U}}(k)$ for $k = 0, 1, 2$, which are obtained by our proposed control scheme.}
\label{fig3}
\end{figure}

\begin{figure}[htbp]
\centering
\includegraphics[scale=0.3]{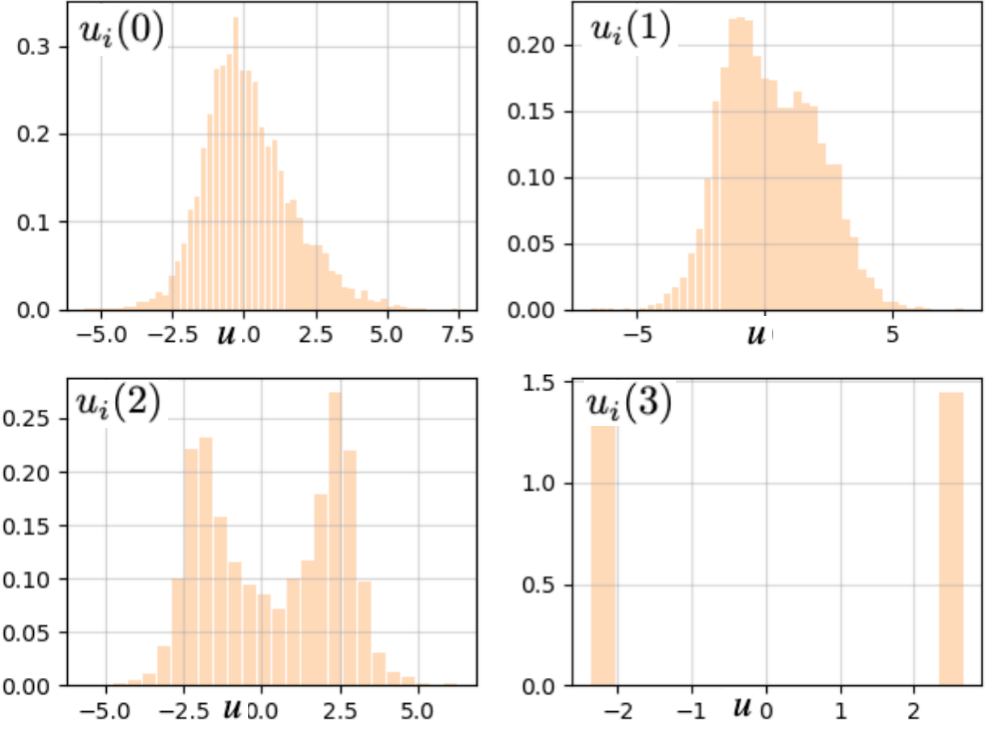}
\centering
\caption{The histograms of $u_{i}(k)$ at time step $k$ for each agent $i$. The upper left and right figures are $u_{i}(0)$ and $u_{i}(1), i = 1, \cdots, 5000$ respectively. The lower left and right figures are $u_{i}(2)$ and $u_{i}(3)$ respectively.}
\label{fig4}
\end{figure}

\begin{figure}[htbp]
\centering
\includegraphics[scale=0.40]{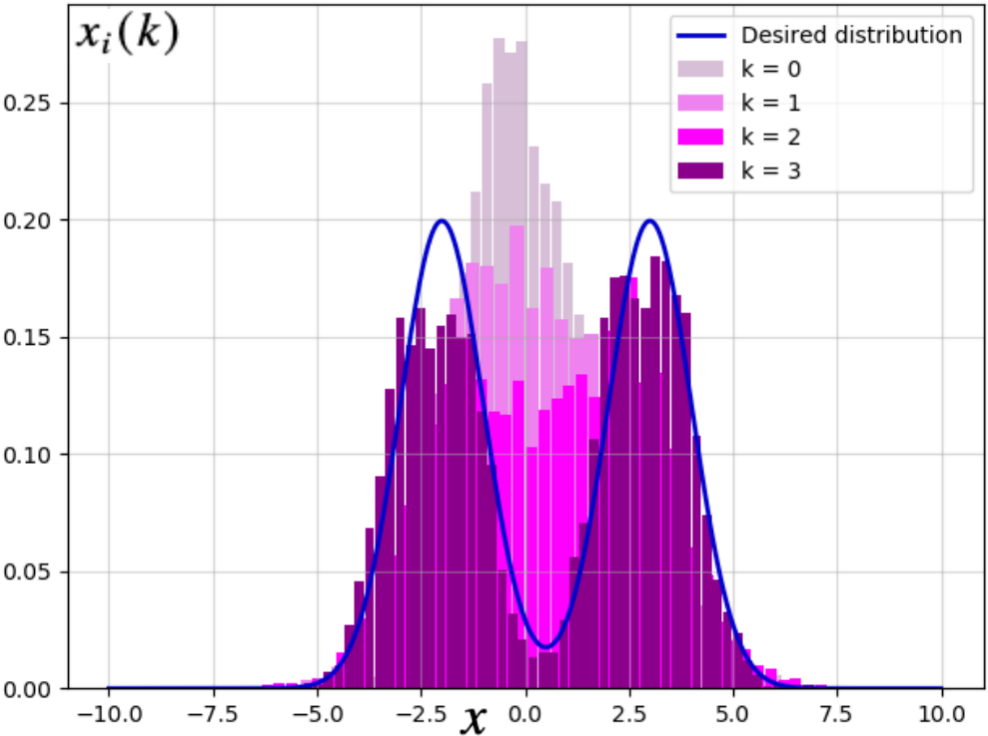}
\centering
\caption{The histogram of the terminal distribution of $x_{i}(k), i = 1, \ldots, 5000$ compared to the desired distribution. It is close to the specified terminal distribution \eqref{qt2}.}
\label{fig5}
\end{figure}

The system parameter $a(k)$ follows the Laplace distribution $\mathcal{C}\left( 0.5, 0.1 \right)$ for all $k=0,1,2,3$. 
As a result, for some agent dynamics $|a_i(k)|>1$, which could cause their corresponding systems to be unstable.
The results are given in Figures~\ref{fig1}-\ref{fig5}. 
The optimal values of $c(k) = 0$ for $k = 0, 1, 2$, and $c(3) = 0.16$.
Since $c(3) \neq 0$, the corresponding $\tilde{\mathscr{U}}(3)$ 
is not positive definite, with its determinant being zero.
This causes the control input $u(3)$ to be a discrete distribution supported on two points $-2.38, 2.69$, with corresponding probability values $0.472, 0.528$, respectively.
From Figure~\ref{fig5}, we note that the terminal discrete 
distribution using the proposed algorithm is close to the desired continuous one. 
The total processing time is $8.528$ seconds for determining the control inputs of all agents across all time steps (tested on an Apple M4 Max chip), which demonstrates the computational efficiency of the algorithm for controlling large-population ensembles.

\section{Conclusions} \label{sec6}

In this paper, we have addressed the control of a large group of agents through a generalized distribution steering problem. 
We present a novel formulation for the dynamics of the ensemble system. 
Unlike traditional approaches, this distribution steering problem does not rely on assuming Gaussian distributions for the system states or the parameters. 
Instead, we only require the existence of finite power moments up to order $2n$, where $n$ is some positive integer. To our knowledge, this algorithm is the first viable solution for the distribution steering of a large swarm robot group governed by linear ODE systems with uncertainty in the system parameter.

Looking forward, future research will aim to extend the results of this paper to multiple-dimensional and/or nonlinear systems. 
However, such an extension poses major challenges. 
Ensuring the existence of the control inputs, namely, control inputs with non-negative probability density functions, poses a significant hurdle necessitating the application 
of tools from real algebraic geometry, e.g., Positivstellensatz.

\bibliographystyle{ieeetr}
\bibliography{refs}

@article{li2010ensemble,
  title={Ensemble control of finite-dimensional time-varying linear systems},
  author={Li, Jr-Shin},
  journal={IEEE Transactions on Automatic Control},
  volume={56},
  number={2},
  pages={345--357},
  year={2010},
  publisher={IEEE}
}

@article{wu2022group,
  title={Group Steering: Approaches Based on Power Moments},
  author={Wu, Guangyu and Lindquist, Anders},
  journal={arXiv preprint arXiv:2211.13370},
  year={2022}
}

@article{wu2023general,
  title={General Distribution Steering: A Sub-Optimal Solution by Convex Optimization},
  author={Wu, Guangyu and Lindquist, Anders},
  journal={arXiv preprint arXiv:2301.06227},
  year={2023}
}

@book{schmudgen2017moment,
  title={{The Moment Problem}},
  author={Schm{\"u}dgen, Konrad},
  publisher={Graduate Texts in Mathematics, Springer},
  volume={277},
  year={2017}
}

@book{boyd2004convex,
  title={{Convex Optimization}},
  author={Boyd, Stephen and Vandenberghe, Lieven},
  year={2004},
  publisher={Cambridge University Press}
}

@article{wu2023non,
  title={Non-{G}aussian {B}ayesian filtering by density parametrization using power moments},
  author={Wu, Guangyu and Lindquist, Anders},
  journal={Automatica},
  volume={153},
  pages={111061},
  year={2023},
  publisher={Elsevier}
}

@article{biswal2021decentralized,
  title={Decentralized Control of Multiagent Systems Using Local Density Feedback},
  author={Biswal, Shiba and Elamvazhuthi, Karthik and Berman, Spring},
  journal={IEEE Transactions on Automatic Control},
  volume={67},
  number={8},
  pages={3920--3932},
  year={2021},
  publisher={IEEE}
}

@book{bullo2009distributed,
  title={Distributed Control of Robotic Networks: A Mathematical Approach to Motion Coordination Algorithms},
  author={Bullo, Francesco and Cort{\'e}s, Jorge and Martinez, Sonia},
  volume={27},
  year={2009},
  publisher={Princeton University Press}
}

@incollection{mesbahi2010graph,
  title={Graph Theoretic Methods in Multiagent Networks},
  author={Mesbahi, Mehran and Egerstedt, Magnus},
  booktitle={Graph Theoretic Methods in Multiagent Networks},
  year={2010},
  publisher={Princeton University Press}
}

@article{elamvazhuthi2018optimal,
  title={Optimal transport over deterministic discrete-time nonlinear systems using stochastic feedback laws},
  author={Elamvazhuthi, Karthik and Grover, Piyush and Berman, Spring},
  journal={IEEE Control Systems Letters},
  volume={3},
  number={1},
  pages={168--173},
  year={2018},
  publisher={IEEE}
}

@article{li2009ensemble,
  title={{Ensemble control of Bloch equations}},
  author={Li, Jr-Shin and Khaneja, Navin},
  journal={IEEE Transactions on Automatic Control},
  volume={54},
  number={3},
  pages={528--536},
  year={2009},
  publisher={IEEE}
}

@article{bakolas2019dynamic,
  title={{Dynamic output feedback control of the Liouville equation for discrete-time SISO linear systems}},
  author={Bakolas, Efstathios},
  journal={IEEE Transactions on Automatic Control},
  volume={64},
  number={10},
  pages={4268--4275},
  year={2019},
  publisher={IEEE}
}

@article{okamoto2019optimal,
  title={Optimal stochastic vehicle path planning using covariance steering},
  author={Okamoto, Kazuhide and Tsiotras, Panagiotis},
  journal={IEEE Robotics and Automation Letters},
  volume={4},
  number={3},
  pages={2276--2281},
  year={2019},
  publisher={IEEE}
}

@article{okamoto2018optimal,
  title={Optimal covariance control for stochastic systems under chance constraints},
  author={Okamoto, Kazuhide and Goldshtein, Maxim and Tsiotras, Panagiotis},
  journal={IEEE Control Systems Letters},
  volume={2},
  number={2},
  pages={266--271},
  year={2018},
  publisher={IEEE}
}

@article{balci2020covariance,
  title={{Covariance steering of discrete-time stochastic linear systems based on Wasserstein distance terminal cost}},
  author={Balci, Isin M and Bakolas, Efstathios},
  journal={IEEE Control Systems Letters},
  volume={5},
  number={6},
  pages={2000--2005},
  year={2020},
  publisher={IEEE}
}

@article{liu2022optimal,
  title={Optimal covariance steering for discrete-time linear stochastic systems},
  author={Liu, Fengjiao and Rapakoulias, George and Tsiotras, Panagiotis},
  journal={IEEE Transactions on Automatic Control},
  month = {April},
  year={2025}
}

@article{chen2015optimal,
  title={{Optimal steering of a linear stochastic system to a final probability distribution, Part I}},
  author={Chen, Yongxin and Georgiou, Tryphon T and Pavon, Michele},
  journal={IEEE Transactions on Automatic Control},
  volume={61},
  number={5},
  pages={1158--1169},
  year={2015},
  publisher={IEEE}
}

@article{chen2015optimal2,
  title={{Optimal steering of a linear stochastic system to a final probability distribution, Part II}},
  author={Chen, Yongxin and Georgiou, Tryphon T and Pavon, Michele},
  journal={IEEE Transactions on Automatic Control},
  volume={61},
  number={5},
  pages={1170--1180},
  year={2015},
  publisher={IEEE}
}

@article{yongxin2018optimal,
  title={{Optimal steering of a linear stochastic system to a final probability distribution, Part III}},
  author={Chen, Yongxin and Georgiou, Tryphon T and Pavon, Michele},
  journal={IEEE Transactions on Automatic Control},
  volume={63},
  number={9},
  pages={3112--3118},
  year={2018},
  publisher={Institute of Electrical and Electronics Engineers Inc.}
}

@inproceedings{knaup2023covariance,
    author = {Knaup, Jacob and Tsiotras, Panagiotis},
    title = {Covariance Steering for Systems Subject to Unknown Parameters},
    booktitle = {62nd IEEE Conference on Decision and Control},
address = {Marina Bay Sands, Singapore},
month = {13--15,},
    year = {2023},
pages = {1790--1795}
}

@inproceedings{collins1985covariance,
  author    = {E. Collins and R. Skelton},
  title     = {Covariance control of discrete systems},
  booktitle = {Proc. IEEE Conf. Decision Control},
  address   = {Lauderdale, FL},
  year      = {1985},
  pages     = {542--547},
}

@article{collins1987theory,
  title={A theory of state covariance assignment for discrete systems},
  author={Collins, EMMAN and Skelton, R},
  journal={IEEE Transactions on Automatic Control},
  volume={32},
  number={1},
  pages={35--41},
  year={1987},
  publisher={IEEE}
}

@article{hsieh1990all,
  title={All covariance controllers for linear discrete-time systems},
  author={Hsieh, Chen and Skelton, Robert E},
  journal={IEEE Transactions on Automatic Control},
  volume={35},
  number={8},
  pages={908--915},
  year={1990},
  publisher={IEEE}
}

@article{xu1992improved,
  title={An improved covariance assignment theory for discrete systems},
  author={Xu, J-H and Skelton, Robert E},
  journal={IEEE Transactions on Automatic Control},
  volume={37},
  number={10},
  pages={1588--1591},
  year={1992},
  publisher={IEEE}
}

@article{deshmukh2018mean,
  title={Mean-field stabilization of {Markov} chain models for robotic swarms: Computational approaches and experimental results},
  author={Deshmukh, Vaibhav and Elamvazhuthi, Karthik and Biswal, Shiba and Kakish, Zahi and Berman, Spring},
  journal={IEEE Robotics and Automation Letters},
  volume={3},
  number={3},
  pages={1985--1992},
  year={2018},
  publisher={IEEE}
}

@inproceedings{elamvazhuthi2018mean,
  author = {Elamvazhuthi, Karthik and Biswal, Shiba and Berman, Spring},
  title = {Mean-field stabilization of robotic swarms to probability distributions with disconnected supports},
  booktitle = {Proc. American Control Conference},
  address = {Milwaukee, WI},
  year = {2018},
  pages = {885--892},
}

@inproceedings{rubenstein2013collective,
  title={Collective transport of complex objects by simple robots: theory and experiments},
  author={Rubenstein, Michael and Cabrera, Adrian and Werfel, Justin and Habibi, Golnaz and McLurkin, James and Nagpal, Radhika},
  booktitle={International Conference on Autonomous Agents and 
             Multi-agent Systems},
  pages={47--54},
  year={2013}
}

@inproceedings{ikumapayi2024swarm,
  title={Swarm Robotics in a Sustainable Warehouse Automation: Opportunities, Challenges and Solutions},
  author={Ikumapayi, Omolayo Michael and Laseinde, Opeyeolu Timothy and Elewa, Remilekun R and Ogedengbe, Temitayo Samson and Akinlabi, Esther Titilayo},
  booktitle={{E3S Web of Conferences}},
  volume={552},
  pages={01080},
  year={2024},
  organization={EDP Sciences}
}

@inproceedings{abhang2024swarm,
  title={Swarm Intelligence for Multi-Robot Coordination in Agricultural Automation},
  author={Abhang, LB and Gummadi, Annapurna and Changala, Ravindra and Vuyyuru, Veera Ankalu and Sabareesh, R and Raj, I Infant},
  booktitle={10th International Conference on Advanced Computing and Communication Systems (ICACCS)},
  volume={1},
  pages={455--460},
  year={2024},
}

@article{biswal2020stabilization,
  title={Stabilization of nonlinear discrete-time systems to target measures using stochastic feedback laws},
  author={Biswal, Shiba and Elamvazhuthi, Karthik and Berman, Spring},
  journal={IEEE Transactions on Automatic Control},
  volume={66},
  number={5},
  pages={1957--1972},
  year={2020},
  publisher={IEEE}
}

@inproceedings{nodozi2023physics,
  author = {Nodozi, Iman and O’Leary, Jared and Mesbah, Ali and Halder, Abhishek},
  title = {A physics-informed deep learning approach for minimum effort stochastic control of colloidal self-assembly},
  booktitle = {Proc. American Control Conference},
  address = {San Diego, CA},
  year = {2023},
  pages = {609--615},
}

@article{caluya2021wasserstein,
  title={{Wasserstein proximal algorithms for the Schr{\"o}dinger bridge problem: Density control with nonlinear drift}},
  author={Caluya, Kenneth F and Halder, Abhishek},
  journal={IEEE Transactions on Automatic Control},
  volume={67},
  number={3},
  pages={1163--1178},
  year={2021},
  publisher={IEEE}
}

@inproceedings{sivaramakrishnan2022distribution,
author = {Sivaramakrishnan, Vignesh and Pilipovsky, Joshua and Oishi, Meeko and Tsiotras, Panagiotis},
  title     = {Distribution steering for discrete-time linear systems with general disturbances using characteristic functions},
  booktitle = {Proc. American Control Conference},
  address   = {Atlanta, GA},
  year      = {2022},
  pages     = {4183--4190},
}

@article{georgiou2003kullback,
  title={{Kullback-Leibler approximation of spectral density functions}},
  author={Georgiou, Tryphon T and Lindquist, Anders},
  journal={IEEE Transactions on Information Theory},
  volume={49},
  number={11},
  pages={2910--2917},
  year={2003},
  publisher={IEEE}
}

@inproceedings{byrnes2002identifiability,
  title={Identifiability of shaping filters from covariance lags, cepstral windows and {Markov} parameters},
  author={Byrnes, Christopher I and Enqvist, Per and Lindquist, Anders},
  booktitle={Proceedings of the 41st IEEE Conference on Decision and Control, 2002.},
  volume={1},
  pages={246--251},
  year={2002},
  organization={IEEE}
}

@article{byrnes2001finite,
  title={From finite covariance windows to modeling filters: A convex optimization approach},
  author={Byrnes, Christopher I and Gusev, Sergei V and Lindquist, Anders},
  journal={SIAM Review},
  volume={43},
  number={4},
  pages={645--675},
  year={2001},
  publisher={SIAM}
}

@article{dorigo2021swarm,
  title={Swarm robotics: Past, present, and future [point of view]},
  author={Dorigo, Marco and Theraulaz, Guy and Trianni, Vito},
  journal={Proceedings of the IEEE},
  volume={109},
  number={7},
  pages={1152--1165},
  year={2021},
  publisher={IEEE}
}

@article{dorigo2020reflections,
  title={Reflections on the future of swarm robotics},
  author={Dorigo, Marco and Theraulaz, Guy and Trianni, Vito},
  journal={Science Robotics},
  volume={5},
  number={49},
  pages={eabe4385},
  year={2020},
  publisher={American Association for the Advancement of Science}
}

@article{zhang2020modeling,
  title={Modeling collective behaviors: A moment-based approach},
  author={Zhang, Silun and Ringh, Axel and Hu, Xiaoming and Karlsson, Johan},
  journal={IEEE Transactions on Automatic Control},
  volume={66},
  number={1},
  pages={33--48},
  year={2020},
}

@article{pitman1977occupation,
  title={Occupation measures for {Markov} chains},
  author={Pitman, JW},
  journal={Advances in Applied Probability},
  volume={9},
  number={1},
  pages={69--86},
  year={1977},
  publisher={Cambridge University Press}
}

@article{elamvazhuthi2019mean,
  title={Mean-field models in swarm robotics: A survey},
  author={Elamvazhuthi, Karthik and Berman, Spring},
  journal={Bioinspiration \& Biomimetics},
  volume={15},
  number={1},
  pages={015001},
  year={2019},
  publisher={IOP Publishing}
}

@article{glass2009infinite, 
title={Infinite Dimensional controllability}, 
journal={Encyclopedia of Complexity and Systems Science}, 
author={Glass, Olivier},
year={2009}, 
pages={4804–4820}}

@book{matyas1999generalized,
  title={Generalized Method of Moments Estimation},
  author={M{\'a}ty{\'a}s, L{\'a}szl{\'o}},
  volume={5},
  year={1999},
  publisher={Cambridge University Press}
}

@article{ringh2023mean,
  title={Mean field type control with species dependent dynamics via structured tensor optimization},
  author={Ringh, Axel and Haasler, Isabel and Chen, Yongxin and Karlsson, Johan},
  journal={IEEE Control Systems Letters},
  volume = {7},
pages = {2898--2903},
  year={2023},
}

@inproceedings{wu2023density,
  author={Wu, Guangyu and Lindquist, Anders},
  title={Density steering by power moments},
  booktitle = {Proc. IFAC World Congress},
  address   = {Yokohama, Japan},
  year      = {2023},
  pages     = {3423--3428},
}

@article{liu2024optimal,
  title={Optimal covariance steering for continuous-time linear stochastic systems with multiplicative noise},
  author={Liu, Fengjiao and Tsiotras, Panagiotis},
  journal={IEEE Transactions on Automatic Control},
  volume={69},
  number={10},
  pages={7247--7254},
  year={2024},
  publisher={IEEE}
}

@article{yu2024optimal,
  title={Optimal Covariance Steering of Linear Stochastic Systems with Hybrid Transitions},
  author={Yu, Hongzhe and Franco, Diana Frias and Johnson, Aaron M and Chen, Yongxin},
  journal={arXiv preprint arXiv:2410.13222},
  year={2024}
}

\end{document}